\documentclass{article}
\usepackage{amssymb}
\usepackage{graphicx}
\usepackage{amsmath}

\input{tcilatex}
\begin{document}

\begin{center}
\textbf{\ \ Pseudo-differential operators in vector-valued spaces and
applications }

\textbf{VELI SHAKHMUROV}\ \ \ \ \ \ \ \ 

\bigskip Okan University, Department of Mechanical Engineering, Akfirat,
Tuzla 34959 Istanbul, Turkey, E-mail: veli.sahmurov@okan.edu.tr;

\ \ \ \ \ \ \ \ \ \ \ \ \ \ 
\end{center}

\begin{quote}
\ \ \ \ \ \ \ \ \ \ \ \ \ \ \ \ \ \ \ \ \ \ \ \ \ \ \ \ \ \ \ \ \ \ \ \ \ \
\ \ \ \ \ \ \ \ \ \ \ \ \ \ \ 
\end{quote}

\begin{center}
\textbf{ABSTRACT}
\end{center}

\begin{quote}
\ \ \ \ \ \ \ \ \ \ \ \ \ \ \ \ \ 
\end{quote}

\ \ \ Pseudo-differential operator equations with parameter are studied.
Uniform $L_{p}$-separability properties and resolvent estimates are obtained
in terms of fractional derivatives. Moreover, maximal regularity properties
of the pseudo-differential abstract parabolic equation are established.
Particularly, it is proven that the operators generated by these
pseudo-differential equations are positive and aso are generators of
analytic semigroups. As an application, the anisotropic parameter dependent
pseudo-differential equations and the system of pseudo-differential
equations are studied.\ \ \ \ \ 

\begin{center}
\bigskip\ \ \textbf{AMS: 47GXX, 35JXX, 47FXX, 47DXX, 43AXX}
\end{center}


\begin{center}
\textbf{1. Introduction, notations and background }
\end{center}

Differential-operator equations (DOEs) have found many applications in PDEs
and pseudo-differential equations (PsDEs) (see e.g. $\left[ 1\right] ,$ $%
\left[ 5\right] ,$ $\left[ 11\right] ,$ $\left[ 12\right] ,$ $\left[ 16-19%
\right] ,$ $\left[ 24\right] $ ). Regularity properties of PsDEs have been
studied extensively by many researchers; see e.g. $\left[ \text{6, 10}\right]
,$ $\left[ \text{21-22}\right] \ $and the references therein. The
boundedness of PsDEs in Sobolev spaces have been treated e.g. in $\left[ 10%
\right] ,$ $\left[ 14\right] ,$ $\left[ 22\right] .$ Moreover, the
smoothness of PsDEs with bounded operator coefficients have been explored
e.g. in $\left[ 8\right] $, $\left[ 15\right] .$ In contrast to $\left[ 8%
\right] $, $\left[ 15\right] ,$ the PsDE considered here contain unbounded
operators and parameters. In particular, the main objective of the present
paper is to discuss the uniform $L_{p}\left( R^{n};E\right) -$maximal
regularity of elliptic pseudo-differential operator equations (PsDOEs) with
parameters%
\begin{equation}
P_{t}\left( D\right) u+Au+\dsum\limits_{\left\vert \alpha \right\vert
<m}t\left( \alpha \right) A_{\alpha }\left( x\right) D^{\alpha }u+\lambda
u=f\left( x\right) ,\text{ }x\in R^{n},  \tag{1.1}
\end{equation}%
where $P_{t}\left( D\right) $ is the pseudo-differential operator, $A$ and $%
A_{\alpha }\left( x\right) $ are linear operators in a Banach space $E$,
for\ $\alpha _{i}\in \left[ 0,\infty \right) $, $\alpha =\left( \alpha
_{1},\alpha _{2},...\alpha _{n}\right) $ and $D^{\alpha }=D_{1}^{\alpha
_{1}}D_{2}^{\alpha _{2}}...D_{n}^{\alpha _{n}}$ are the Liouville
derivatives; $m$ is a positive number, $t_{k}$ are positive, $\lambda $ is a
complex parameter, $t=\left( t_{1},t_{2},...,t_{n}\right) $ and $t\left(
\alpha \right) =\dprod\limits_{k=1}^{n}t_{k}^{\frac{\alpha _{k}}{m}};$ $%
L_{p}\left( \Omega ;E\right) $ denotes the space of strongly measurable $E$%
-valued functions that are defined on the measurable subset $\Omega \subset
R^{n}$ with the norm given by

\begin{equation*}
\left\Vert f\right\Vert _{L_{p}\left( \Omega ;E\right) }=\left(
\int\limits_{\Omega }\left\Vert f\left( x\right) \right\Vert
_{E}^{p}dx\right) ^{\frac{1}{p}},\text{ }1\leq p<\infty \ .
\end{equation*}%
We prove that problem $\left( 1.1\right) $ has a maximal regular unique
solution and the following uniform coercive estimate holds 
\begin{equation}
\dsum\limits_{\left\vert \alpha \right\vert \leq m}t\left( \alpha \right)
\left\vert \lambda \right\vert ^{1-\frac{\left\vert \alpha \right\vert }{m}%
}\left\Vert D^{\alpha }u\right\Vert _{L_{p}\left( R^{n};E\right)
}+\left\Vert Au\right\Vert _{L_{p}\left( R^{n};E\right) }\leq C\left\Vert
f\right\Vert _{L_{p}\left( R^{n};E\right) }  \tag{1.2}
\end{equation}%
for $f\in L_{p}\left( R^{n};E\right) $, $\lambda \in S_{\varphi },$ where $%
S_{\varphi }$ is a set of complex numbers that is related with the spectrum
of the operator $A.$\ The estimate $\left( 1.2\right) $ implies that the
operator $O_{t}$ generated by $\left( 1.1\right) $ has a bounded inverse
from $L_{p}\left( R^{n};E\right) $ into the space $H_{p}^{m}\left(
R^{n};E\left( A\right) ,E\right) ,$ which will be defined subsequently.
Particularly, from the estimate $\left( 1.2\right) $ we obtain that the
operator $O_{t}$ is uniformly positive in $L_{p}\left( R^{n};E\right) .$ By
using this property we prove the uniform well posedness of the Cauchy
problem for the following parabolic PsDOE with parameter

\begin{equation}
\frac{\partial u}{\partial y}+P_{t}\left( D\right) u+Au=f\left( y,x\right) ,%
\text{ }u(0,x)=0,  \tag{1.3}
\end{equation}%
in $E$-valued mixed spaces $L_{\mathbf{p}},$ $\mathbf{p=}\left(
p,p_{1}\right) $. In other words, we show that problem $\left( 1.3\right) $
has a unique solution 
\begin{equation*}
u\in W_{\mathbf{p}}^{1,m}\left( R_{+}^{n+1};E\left( A\right) ,E\right)
\end{equation*}%
for $f\in L_{\mathbf{p}}\left( R_{+}^{n+1};E\right) $\ satisfying the
following uniform coercive estimate 
\begin{equation}
\left\Vert \frac{\partial u}{\partial y}\right\Vert _{L_{_{\mathbf{p}%
}}\left( R_{+}^{n+1};E\right) }+\left\Vert P_{t}\left( D\right) u\right\Vert
_{L_{\mathbf{p}}\left( R_{+}^{n+1};E\right) }+  \tag{1.4}
\end{equation}

\begin{equation*}
\left\Vert Au\right\Vert _{L_{\mathbf{p}}\left( R_{+}^{n+1};E\right) }\leq
M\left\Vert f\right\Vert _{L_{\mathbf{p}}\left( R_{+}^{n+1};E\right) }.
\end{equation*}%
Note that, constants $C$, $M$ in estimates $\left( 1.2\right) $ and $\left(
1.4\right) $ are independent of parameters. As an application in this paper
the following are established: (a) maximal regularity properties of the
anisotropic elliptic PsDE in mixed $L_{\mathbf{p}},$ $\mathbf{p=}\left(
p_{1},p\right) $ spaces; (b) coercive properties of the system of PsDEs of
infinite order in $L_{p}$ spaces.

Let $\left( \Omega ;\Sigma ,\mu \right) $ be a complete probability space, $%
1\leq p<\infty .$ $L_{p}\left( \Omega ;\Sigma ,\mu ,E\right) $ denotes $\mu
- $measurable $E-$valued Bochner space with norm 
\begin{equation*}
\left\Vert f\right\Vert _{L_{p}\left( \Omega ;\Sigma ,\mu ,E\right) }=\left(
\int\limits_{\Omega }\left\Vert f\left( x\right) \right\Vert _{E}^{p}d\mu
\right) ^{\frac{1}{p}}.
\end{equation*}

A Banach space $E$ is called UMD space (see $\left[ \text{7, \S\ 5}\right] $%
) if $E$-valued martingale difference sequences are unconditional in $%
L_{p}\left( \Omega ;\Sigma ,\mu ,E\right) $ for $p\in \left( 1,\infty
\right) $, i.e., there exists a positive constant $C_{p}$ such that for any
martingale $\left\{ f_{k}\right\} $ any choice of signs $\left\{ \varepsilon
_{k}\right\} \in \left\{ -1,1\right\} $, $k\in \mathbb{N}$ and $N\in \mathbb{%
N}$%
\begin{equation*}
\left\Vert f_{0}+\sum\limits_{k=1}^{N}\varepsilon _{k}\left(
f_{k}-f_{k-1}\right) \right\Vert _{L_{p}\left( \Omega ,\Sigma ,\mu ,E\right)
}\leq C_{p}\left\Vert f_{N}\right\Vert _{L_{p}\left( \Omega ,\Sigma ,\mu
,E\right) }.
\end{equation*}

It is shown (see $\left[ 3-4\right] $) that\ the Hilbert operator 
\begin{equation*}
\left( Hf\right) \left( x\right) =\lim\limits_{\varepsilon \rightarrow
0}\int\limits_{\left\vert x-y\right\vert >\varepsilon }\frac{f\left(
y\right) }{x-y}dy
\end{equation*}%
is bounded in the space $L_{p}\left( R,E\right) ,$ $p\in \left( 1,\infty
\right) $ for those and only those spaces $E$ which possess the property of
UMD spaces. UMD spaces include e.g. $L_{p}$, $l_{p}$ and Lorentz spaces $%
L_{pq},$ $p,$ $q\in \left( 1,\infty \right) $.

Let $\mathbb{C}$ denote the set of complex numbers and 
\begin{equation*}
S_{\varphi }=\left\{ \lambda ;\text{ \ }\lambda \in \mathbb{C}\text{, }%
\left\vert \arg \lambda \right\vert \leq \varphi \right\} \cup \left\{
0\right\} ,\text{ }0\leq \varphi <\pi .\ 
\end{equation*}

A linear operator\ $A$ is said to be $\varphi $-positive (or positive ) in a
Banach\ space $E$ if $D\left( A\right) $ is dense on $E$ and%
\begin{equation*}
\left\Vert \left( A+\lambda I\right) ^{-1}\right\Vert _{B\left( E\right)
}\leq M\left( 1+\left\vert \lambda \right\vert \right) ^{-1}
\end{equation*}%
for any $\lambda \in S_{\varphi },$ where $\varphi \in \left[ 0\right.
,\left. \pi \right) $, $I$ is the identity operator in $E,$ $B\left(
E\right) $ is the space of bounded linear operators in $E.$ Sometimes $%
A+\lambda I$\ will be written $A+\lambda $ and will be denoted by $%
A_{\lambda }.$ It is known $\left[ \text{20, \S 1.15.1}\right] $ that the
powers\ $A^{\theta }$, $\theta \in \left( -\infty ,\infty \right) $ for a
positive operator $A$ exist$.$

The operator $A\left( h\right) ,$ $h\in Q\subset \mathbb{C}$ is said to be $%
\varphi $-positive (or positive) in\ $E$ uniformly with respect to $h\in Q$
if $D\left( A\left( h\right) \right) $ is independent of $h$, $D\left(
A\left( h\right) \right) $ is dense in $E$ and $\left\Vert \left( A\left(
h\right) +\lambda \right) ^{-1}\right\Vert \leq M\left( 1+\left\vert \lambda
\right\vert \right) ^{-1}$ for all $\lambda \in S_{\varphi },$ $0\leq
\varphi <\pi $, where $M$ does not depend on $h$ and $\lambda .$ Let $%
E\left( A^{\theta }\right) $ denote the space $D\left( A^{\theta }\right) $
with the norm 
\begin{equation*}
\left\Vert u\right\Vert _{E\left( A^{\theta }\right) }=\left( \left\Vert
u\right\Vert ^{p}+\left\Vert A^{\theta }u\right\Vert ^{p}\right) ^{\frac{1}{p%
}},\text{ }1\leq p<\infty ,\text{ }0<\theta <\infty .
\end{equation*}

A set $W\subset B\left( E_{1},E_{2}\right) $ is called $R$-bounded (see e.g. 
$\left[ 23\right] $) if there is a constant $C>0$ such that for all $%
T_{1},T_{2},...,T_{m}\in W$ and $u_{1,}u_{2},...,u_{m}\in E_{1},$ $m\in 
\mathbb{N}$ 
\begin{equation*}
\int\limits_{0}^{1}\left\Vert \sum\limits_{j=1}^{m}r_{j}\left( y\right)
T_{j}u_{j}\right\Vert _{E_{2}}dy\leq C\int\limits_{0}^{1}\left\Vert
\sum\limits_{j=1}^{m}r_{j}\left( y\right) u_{j}\right\Vert _{E_{1}}dy,
\end{equation*}%
where $\left\{ r_{j}\right\} $ is an arbitrary sequence of independent
symmetric $\left\{ -1,1\right\} $-valued random variables on $\left[ 0,1%
\right] $.

The smallest $C$ for which the above estimate holds is called a $R$-bound of
the collection $W$ and is denoted by $R\left( W\right) .$

A set of operators $G_{h}\subset B\left( E_{1},E_{2}\right) $ depending on
parameter $h\in Q\subset \mathbb{C}$ is called uniformly $R$-bounded with
respect to $h$ if there is a constant $C$ independent of $h\in Q$ such that 
\begin{equation*}
\int\limits_{\Omega }\left\Vert \sum\limits_{j=1}^{m}r_{j}\left( y\right)
T_{j}\left( h\right) u_{j}\right\Vert _{E_{2}}dy\leq C\int\limits_{\Omega
}\left\Vert \sum\limits_{j=1}^{m}r_{j}\left( y\right) u_{j}\right\Vert
_{E_{1}}dy
\end{equation*}%
for all $T_{1}\left( h\right) ,T_{2}\left( h\right) ,...,T_{m}\left(
h\right) \in G_{h}$ and $u_{1,}u_{2},...,u_{m}\in E_{1},$ $m\in \mathbb{N}$.

It implies that 
\begin{equation*}
\sup\limits_{h\in Q}R\left( G_{h}\right) \leq C.
\end{equation*}

\ \ The operator $A$ is said to be $R$-positive in a Banach space $E$ if the
set $\left\{ \lambda \left( A+\lambda \right) ^{-1}\text{: }\lambda \in
S_{\varphi }\right\} $ is $R$-bounded.

A positive operator $A\left( h\right) $ is said to be uniformly $R$-positive
in a Banach space $E$ if there exists $\varphi \in \left[ 0\right. ,\left.
\pi \right) $ such that the set%
\begin{equation*}
\left\{ \lambda \left( A\left( h\right) +\lambda \right) ^{-1}:\lambda \in
S_{\varphi }\right\}
\end{equation*}%
is uniformly $R$-bounded. Let $S\left( R^{n};E\right) $ denote the $E-$%
valued Schwartz class, i.e., the space of all $E$-valued rapidly decreasing
smooth functions on $R^{n}$ equipped with its usual topology generated by
seminorms. For $E=\mathbb{C}$ this space will be denoted by $S=S\left(
R^{n}\right) $. $S^{\prime }\left( E\right) =S^{\prime }\left(
R^{n};E\right) $ denotes the space of linear continuous mappings from $S$
into\ $E$ and is called $E$-valued Schwartz distributions. For any $\alpha
=\left( \alpha _{1},\alpha _{2},...\alpha _{n}\right) $, \ $\alpha _{i}\in %
\left[ 0,\infty \right) $ the function $\left( i\xi \right) ^{\alpha }$ will
be defined such that

{\large 
\begin{equation*}
\left( i\xi \right) ^{\alpha }=\left\{ 
\begin{array}{c}
\left( i\xi _{1}\right) ^{\alpha _{1}}...\left( i\xi _{n}\right) ^{\alpha
_{n}}\text{, \ \ }\xi _{1}\xi _{2},...,\xi _{n}\neq 0 \\ 
0\text{,\ \ \ \ \ }\xi _{1},\xi _{2},...,\xi _{n}=0,%
\end{array}%
\right.
\end{equation*}%
}

where 
\begin{equation*}
\left( i\xi _{k}\right) ^{\alpha _{k}}=\exp \left[ \alpha _{k}\left( \ln
\left\vert \xi _{k}\right\vert +i\pi \text{ sgn }\xi _{k}/2\right) \right] 
\text{, }k=1,2,...,n.
\end{equation*}%
The Liouville derivatives $D^{\alpha }u$ of an $E$-valued function $u$ are
defined similarly to the case of scalar functions $\left[ \text{13}\right] .$

$C\left( \Omega ;E\right) $ and $C^{\left( m\right) }\left( \Omega ;E\right) 
$\ will denote the spaces of $E-$valued bounded uniformly strongly
continuous and $m$ times continuously differentiable functions on $\Omega $,
respectively. Let $F$ and $F^{-1}$ denote the Fourier and inverse Fourier
transforms defined as%
\begin{equation*}
Fu=\left( 2\pi \right) ^{-\frac{n}{2}}\dint\limits_{R^{n}}\left[ \exp \left(
x,\xi \right) \right] u\left( x\right) dx,\text{ }F^{-1}u=\left( 2\pi
\right) ^{-\frac{n}{2}}\dint\limits_{R^{n}}\left[ \exp -\left( x,\xi \right) %
\right] u\left( \xi \right) d\xi ,
\end{equation*}%
where 
\begin{equation*}
x=\left( x_{1},x_{2},...,x_{n}\right) ,\text{ }\xi =\left( \xi _{1},\xi
_{2},...,\xi _{n}\right) \in R^{n}\text{, }\left( x,\xi \right)
=\dsum\limits_{k=1}^{n}x_{k}\xi _{k}.
\end{equation*}%
Through this section, the Fourier transformation of a function $u$ will be
denoted by $\hat{u}.$ It is known that 
\begin{equation*}
\ F\left( D_{x}^{\alpha }u\right) =\left( i\xi _{1}\right) ^{\alpha
_{1}}...\left( i\xi _{n}\right) ^{\alpha _{n}}\hat{u},\ D_{\xi }^{\alpha
}\left( F\left( u\right) \right) =F\left[ \left( -ix_{n}\right) ^{\alpha
_{1}}...\left( -ix_{n}\right) ^{\alpha _{n}}u\ \ \right] \ \ \ \ \ 
\end{equation*}%
for all $u\in S^{\prime }\left( R^{n};E\right) $. Let $E_{1}$ and $E_{2}$ be
two Banach spaces. $B\left( E_{1},E_{2}\right) $ denotes the space of
bounded linear operators from $E_{1}$ to $E_{2}$.\ A function $\Psi \in
C\left( R^{n};B\left( E_{1},E_{2}\right) \right) $ is called a Fourier
multiplier from $L_{p}\left( R^{n};E_{1}\right) $\ to $L_{p}\left(
R^{n};E_{2}\right) $ if the map%
\begin{equation*}
u\rightarrow \Lambda u=F^{-1}\Psi \left( \xi \right) Fu,\text{ }u\in S\left(
R^{n};E_{1}\right)
\end{equation*}%
is well defined and extends to a bounded linear operator 
\begin{equation*}
\Lambda :\ L_{p}\left( R^{n};E_{1}\right) \rightarrow \ L_{p}\left(
R^{n};E_{2}\right) .
\end{equation*}%
The set of all multipliers from\ $L_{p}\left( R^{n};E_{1}\right) $ to\ $%
L_{p}\left( R^{n};E_{2}\right) $\ will be denoted by $M_{p}^{p}\left(
E_{1},E_{2}\right) .$ For $E_{1}=E_{2}=E$\ it is denoted by $M_{p}^{p}\left(
E\right) .$

Let $\Phi _{h}=\left\{ \Psi _{h}\in M_{p}^{p}\left( E_{1},E_{2}\right) ,%
\text{ }h\in Q\right\} $ denote a collection of multipliers depending on the
parameter $h.$

We say that $W_{h}$ is a uniform collection of multipliers if there exists a
positive constant $M$ independent of $h\in Q$ such that

\begin{equation*}
\left\Vert F^{-1}\Psi _{h}Fu\right\Vert _{L_{p}\left( R^{n};E_{2}\right)
}\leq M\left\Vert u\right\Vert _{L_{p}\left( R^{n};E_{1}\right) }\ \ \ \ \ \ 
\end{equation*}%
for all $h\in Q$ and $u\in S\left( R^{n};E_{1}\right) .$

Let $E_{0}$ and $E$ be two Banach spaces and $E_{0}$ be continuously and
densely embedded into $E.$ Let $s\in R$ and $\xi =\left( \xi _{1},\xi
_{2},...,\xi _{n}\right) \in R^{n}.$ Consider the following Liouville-Lions
space

{\large 
\begin{equation*}
H_{p}^{s}(R^{n};E_{0},E)=\left\{ u\right. \ u\in S^{\prime }\left(
R^{n};E_{0}\right) \text{, }F^{-1}\left( 1+\left\vert \xi \right\vert
\right) ^{\frac{s}{2}}Fu\in L_{p}\left( R^{n};E\right) \text{, }
\end{equation*}%
}

\begin{equation*}
\left\Vert u\right\Vert _{H_{p}^{s}\left( R^{n};E_{0},E\right) }=\left\Vert
u\right\Vert _{L_{p}\left( R^{n};E_{0}\right) }+\text{ }\left. \left\Vert
F^{-1}\left( 1+\left\vert \xi \right\vert \right) ^{\frac{s}{2}%
}Fu\right\Vert _{L_{p}\left( R^{n};E\right) }<\infty \right\} .
\end{equation*}

Let $t=\left( t_{1},t_{2},...,t_{n}\right) $ and $t_{k}$ be positive
parameters. We define the following parameterized norm in $%
H_{p}^{s}(R^{n};E_{0},E):$%
\begin{equation*}
\left\Vert u\right\Vert _{H_{p,t}^{s}\left( R^{n};E_{0},E\right)
}=\left\Vert u\right\Vert _{L_{p}\left( R^{n};E_{0}\right) }+\text{ }%
\left\Vert F^{-1}\left[ 1+\left( \dsum\limits_{k=1}^{n}t_{k}^{\frac{2}{s}%
}\xi _{k}^{2}\right) ^{\frac{1}{2}}\right] ^{s}Fu\right\Vert _{L_{p}\left(
R^{n};E\right) }<\infty .
\end{equation*}

Sometimes we use one and the same symbol $C$ without distinction in order to
denote positive constants which may differ from each other even in a single
context. When we want to specify the dependence of such a constant on a
parameter, say $\alpha $, we write $C_{\alpha }$.

By using the techniques of $\left[ \text{9},\text{Theorem 3.7}\right] $ and
reasoning as in $\left[ \text{19, Theorem A}_{0}\right] $ we obtain the
following proposition$.$

\textbf{Proposition A}$_{0}.$ Let $E_{1}$ and $E_{2}$ be two $UMD$ spaces
and 
\begin{equation*}
\Psi _{h}\in C^{n}\left( R^{n}\backslash \left\{ 0\right\} ;B\left(
E_{1},E_{2}\right) \right) .
\end{equation*}
Suppose there is a positive constant $K$ such that

\begin{equation*}
\text{ }\sup\limits_{h\in Q}R\left( \left\{ \left\vert \xi \right\vert
^{\left\vert \beta \right\vert }D^{\beta }\Psi _{h}\left( \xi \right) \text{%
: }\xi \in R^{n}\backslash \left\{ 0\right\} ,\text{ }\beta _{i}\in \left\{
0,1\right\} \right\} \right) \leq K,\text{ }
\end{equation*}

for

\begin{equation*}
\beta =\left( \beta _{1},\beta _{2},...,\beta _{n}\right) ,\text{ }%
\left\vert \beta \right\vert =\dsum\limits_{k=1}^{n}\beta _{k}.
\end{equation*}

Then $\Psi _{h}$ is a uniform collection of multipliers from $L_{p}\left(
R^{n};E_{1}\right) $ to $L_{p}\left( R^{n};E_{2}\right) $ for $p\in \left(
1,\infty \right) .$

\textbf{Proof: }Some steps (Lemma 3.1, Proposition 3.2) of proof $\left[ 
\text{9},\text{Theorem 3.7}\right] $ trivially work for the parameter
dependent case. Other steps (Theorem 3.3, Lemma 3.5) can be easily shown by
replacing 
\begin{equation*}
\left\{ \left\vert \xi \right\vert ^{\left\vert \beta \right\vert }D^{\beta
}\Psi \left( \xi \right) \text{: }\xi \in R^{n}\backslash \left\{ 0\right\}
\right\} \text{ }
\end{equation*}%
with 
\begin{equation*}
\Sigma _{h}=\left\{ \left\vert \xi \right\vert ^{\left\vert \beta
\right\vert }D^{\beta }\Psi _{h}\left( \xi \right) :\xi \in R^{n}\backslash
\left\{ 0\right\} \right\}
\end{equation*}%
\ and by using the uniform $R$-boundedness of the set $\Sigma _{h}$.
However, the parameter dependent analog of Proposition 3.4 in $\left[ \text{9%
}\right] $ is not straightforward. Let $M_{h}$, $M_{h,N}\in
L_{1}^{loc}\left( R^{n},B\left( E_{1},E_{2}\right) \right) $ be Fourier
multipliers from $L_{p}\left( R^{n};E_{1}\right) $ to $L_{p}\left(
R^{n};E_{2}\right) .$ Let $M_{h,N}$ converge to $M_{h}$ in $%
L_{1}^{loc}\left( R^{n},B\left( E_{1},E_{2}\right) \right) $ and $%
T_{h,N}=F^{-1}M_{h,N}F$ be uniformly bounded with respect to $h$ and $N.$
Then the operator $T_{h}=F^{-1}M_{h}F$ is uniformly bounded, so we obtain
the assertion of Proposition A$_{0}.$

The embedding theorems in vector valued spaces play a key role in the theory
of DOEs. For estimating lower order derivatives in terms of interpolation
spaces we use following embedding theorems from $\left[ 17\right] $.

{\large \ }\textbf{Theorem A}$_{1}$\textbf{.} Suppose $E$ is an UMD space, $%
0<t_{k}\leq t_{0}<\infty $, $1<p\leq q<\infty $ and $A$ is an $R$-positive
operator in $E.$ Then for $s$ $\in \left( 0,\infty \right) $ with $\varkappa
=\left\vert \alpha \right\vert +n\left( \frac{1}{p}-\frac{1}{q}\right) \leq
s $, $0\leq \mu \leq 1-\varkappa $ the embedding 
\begin{equation*}
D^{\alpha }H_{p}^{s}\left( R^{n};E\left( A\right) ,E\right) \subset
L_{q}\left( R^{n};E\left( A^{1-\varkappa -\mu }\right) \right)
\end{equation*}%
is continuous and there exists a constant \ $C_{\mu }$ \ $>0$, depending
only on $\mu $ such that 
\begin{equation*}
t\left( \alpha \right) \left\Vert D^{\alpha }u\right\Vert _{L_{q}\left(
R^{n};E\left( A^{1-\varkappa -\mu }\right) \right) }\leq C_{\mu }\left[
h^{\mu }\left\Vert u\right\Vert _{H_{p,t}^{s}\left( R^{n};E\left( A\right)
,E\right) }+h^{-\left( 1-\mu \right) }\left\Vert u\right\Vert _{L_{p}\left(
R^{n};E\right) }\right]
\end{equation*}%
for all $u\in H_{p}^{s}\left( R^{n};E\left( A\right) ,E\right) $ and $%
0<h\leq h_{0}<\infty $ $.$

\begin{center}
2. \textbf{PsDOE with parameters in Banach spaces}
\end{center}

\bigskip Consider the principal part of the problem $\left( 1.1\right) :$ 
\begin{equation}
\left( L_{t}+\lambda \right) u=P_{t}\left( D\right) u+Au+\lambda u=f\left(
x\right) ,\text{ }x\in R^{n},  \tag{2.1}
\end{equation}%
where $P_{t}\left( D\right) $ is the pseudo-differential operator defined by 
\begin{equation}
P_{t}\left( D\right) u=F^{-1}P_{t}\left( \xi \right) \hat{u}\left( \xi
\right) =\left( 2\pi \right) ^{-\frac{n}{2}}\dint\limits_{R^{n}}e^{i\left(
x,\xi \right) }P_{t}\left( \xi \right) \hat{u}\left( \xi \right) d\xi , 
\tag{2.2}
\end{equation}

\textbf{Condition 2.1. }Assume $P_{t}\left( \xi \right) \in $ $S^{m}$ for
some positive number $m,$ i.e.,%
\begin{equation*}
\left\vert D_{\xi }^{\alpha }P_{t}\left( \xi \right) \right\vert \leq
C_{\alpha }\left[ 1+\left( \dsum\limits_{k=1}^{n}t_{k}^{\frac{2}{%
m-\left\vert \alpha \right\vert }}\xi _{k}^{2}\right) ^{\frac{1}{2}}\right]
^{m-\left\vert \alpha \right\vert }
\end{equation*}
for all $\xi \in R^{n}$ and $t_{k}\in \left( 0,\left. t_{0}\right] .\right. $%
\ Suppose $P_{t}\left( \xi \right) \in S_{\varphi _{1}}$ for all $\xi \in
R^{n},$ $t_{k}\in \left( 0,\left. t_{0}\right] ,\right. $ $0\leq \varphi
_{1}<\pi $ and there is a constant $\gamma >0$ such that $\left\vert
P_{t}\left( \xi \right) \right\vert \geq \gamma
\dsum\limits_{k=1}^{n}t_{k}\left\vert \xi _{k}\right\vert ^{m}$.

Let 
\begin{equation*}
Y=H_{p}^{m}\left( R^{n};E\left( A\right) ,E\right) .
\end{equation*}%
\ In this section we prove the following

\textbf{Theorem 2.1.} Assume the Condition 2.1 hold. Suppose $E$ is an UMD
space, $p\in \left( 1,\infty \right) $ and $A$ is an $R$-positive operator
in $E$ with respect to $\varphi \in \left( 0\right. ,\left. \pi \right] .$
Then for $f\in L_{p}\left( R^{n};E\right) $, $\lambda \in S_{\varphi _{2}},$ 
$0\leq \varphi _{2}<\pi -\varphi _{1}$ and $\varphi _{1}+\varphi _{2}\leq
\varphi $ there is a unique solution $u$ of the equation $\left( 2.1\right) $
belonging to $Y$ and the following coercive uniform estimate holds 
\begin{equation}
\dsum\limits_{\left\vert \alpha \right\vert \leq m}t\left( \alpha \right)
\left\vert \lambda \right\vert ^{1-\frac{\left\vert \alpha \right\vert }{m}%
}\left\Vert D^{\alpha }u\right\Vert _{L_{p}\left( R^{n};E\right)
}+\left\Vert Au\right\Vert _{L_{p}\left( R^{n};E\right) }\leq C\left\Vert
f\right\Vert _{L_{p}\left( R^{n};E\right) }.  \tag{2.3}
\end{equation}

\textbf{Proof. }By applying the Fourier transform to equation $\left(
2.1\right) $ we\ obtain 
\begin{equation}
\left[ P_{t}\left( \xi \right) +A+\lambda \right] \hat{u}\left( \xi \right) =%
\hat{f}\left( \xi \right) .  \tag{2.4}
\end{equation}

By construction $\lambda +P_{t}\left( \xi \right) \in S_{\varphi },$ for
all\ $t_{k}\in \left( 0\right. ,\left. t_{0}\right] ,$ $\xi \in R^{n}$ and
the operator\ $A+\lambda +P_{t}\left( \xi \right) $ is invertible in $E$.
So, from $\left( 2.4\right) $ we obtain that the solution of equation $%
\left( 2.1\right) $ can be represented in the form 
\begin{equation}
u\left( x\right) =F^{-1}\left[ A+\lambda +P_{t}\left( \xi \right) \right]
^{-1}\hat{f}.  \tag{2.5}
\end{equation}

By definition of the positive operator $A,$ the inverse of $A^{-1}$ is
bounded in $E$. Then the operator $A$ is a closed linear operator (as an
inverse of bounded linear operator $A^{-1}$). By the differential properties
of the Fourier transform and by using $\left( 2.5\right) $ we have

{\large 
\begin{equation*}
\left\Vert Au\right\Vert _{X}=\left\Vert F^{-1}A\left[ A+\lambda
+P_{t}\left( \xi \right) \right] ^{-1}\hat{f}\right\Vert _{X},
\end{equation*}%
}

\begin{equation*}
\left\Vert D^{\alpha }u\right\Vert _{X}=\left\Vert F^{-1}\xi ^{\alpha }\left[
A+\lambda +P_{t}\left( \xi \right) \right] ^{-1}\hat{f}\right\Vert _{X},
\end{equation*}%
where $X=L_{p}\left( R^{n};E\right) .$\ Hence, it suffices to show that
operator-functions 
\begin{equation*}
\sigma \left( t,\lambda ,\xi \right) =A\left[ A+\lambda +P_{t}\left( \xi
\right) \right] ^{-1}\text{, }
\end{equation*}%
\begin{equation*}
\sigma _{\alpha }\left( t,\lambda ,\xi \right) =t\left( \alpha \right)
\left\vert \lambda \right\vert ^{1-\frac{\left\vert \alpha \right\vert }{m}%
}\xi ^{\alpha }\left[ A+\lambda +P_{t}\left( \xi \right) \right] ^{-1}
\end{equation*}%
are collections of multipliers in $X\ $uniformly with respect to $t_{k}\in
\left( 0,\left. t_{0}\right] \right. $ and $\lambda \in S_{\varphi _{2}}.$
By virtue of $\left[ \text{5, Lemma 2.3}\right] ,$ for $\lambda \in
S_{\varphi _{1}}$ and $\nu \in S_{\varphi _{2}}$ with $\varphi _{1}+\varphi
_{2}<\pi $ there is a positive constant $C$ such that 
\begin{equation}
\left\vert \lambda +\nu \right\vert \geq C\left( \left\vert \lambda
\right\vert +\left\vert \nu \right\vert \right) .  \tag{2.6}
\end{equation}%
By using the positivity properties of operator $A$ we get that 
\begin{equation*}
B\left( \lambda ,t\right) =\left[ A+\lambda +P_{t}\left( \xi \right) \right]
^{-1}
\end{equation*}%
is bounded for all $\xi \in R^{n}$, $\lambda \in S_{\varphi _{1}},$ $%
t_{k}\in \left( 0\right. ,\left. t_{0}\right] $ and 
\begin{equation*}
\left\Vert B\left( \lambda ,t\right) \right\Vert \leq C\left( 1+\left\vert
\lambda +P_{t}\left( \xi \right) \right\vert \right) ^{-1}.
\end{equation*}%
By using Condition 2.1\ and estimate $\left( 2.6\right) $\ we obtain that 
\begin{equation}
\left\Vert B\left( \lambda ,t\right) \right\Vert \leq C\left( 1+\left\vert
\lambda \right\vert +\left\vert P_{t}\left( \xi \right) \right\vert \right)
^{-1}\leq  \tag{2.7}
\end{equation}%
\begin{equation*}
C_{2}\left[ 1+\left\vert \lambda \right\vert
+\dsum\limits_{k=1}^{n}t_{k}\left\vert \xi _{k}\right\vert ^{m}\right] ^{-1}.
\end{equation*}%
Then by using the resolvent properties of positive operators and uniform
estimate $\left( 2.7\right) $\ we obtain 
\begin{equation*}
\left\Vert \sigma \left( t,\lambda ,\xi \right) \right\Vert \leq \left\Vert
I+\left( \lambda +P_{t}\left( \xi \right) \right) \left[ A+\lambda
+P_{t}\left( \xi \right) \right] ^{-1}\right\Vert \leq
\end{equation*}%
\begin{equation*}
1+\left( \left\vert \lambda \right\vert +\left\vert P_{t}\left( \xi \right)
\right\vert \right) \left( 1+\left\vert \lambda \right\vert +\left\vert
P_{t}\left( \xi \right) \right\vert \right) ^{-1}\leq C_{3},
\end{equation*}

where $I$ is an identity operator in $E.$ Moreover,\ by using the well known
inequality 
\begin{equation*}
y_{1}^{\beta _{1}}y_{2}^{\beta _{2}}...y_{n}^{\beta _{n}}\leq C\left(
1+\dsum\limits_{k=1}^{n}y_{k}^{m}\right)
\end{equation*}%
for $\left\vert \beta \right\vert \leq m,$ $y_{k}>0$ and $\beta =\left(
\beta _{1},\beta _{2},...,\beta _{n}\right) $ for all $u\in E$ we have%
\begin{equation*}
\left\Vert \sigma _{\alpha }\left( t,\lambda ,\xi \right) u\right\Vert
_{E}\leq t\left( \alpha \right) \left\vert \lambda \right\vert ^{1-\frac{%
\left\vert \alpha \right\vert }{m}}\left\vert \xi ^{\alpha }\right\vert
\left\Vert B\left( \lambda ,t\right) u\right\Vert _{E}\leq
\end{equation*}%
\begin{equation*}
\left\vert \lambda \right\vert \dprod\limits_{k=1}^{n}\left( t_{k}^{\frac{1}{%
m}}\left\vert \lambda \right\vert ^{\frac{1}{m}}\left\vert \xi
_{k}\right\vert \right) ^{\alpha _{k}}\left\Vert B\left( \lambda ,t\right)
u\right\Vert _{E}\leq C_{\alpha }\left( \left\vert \lambda \right\vert
+\dsum\limits_{k=1}^{n}t_{k}\left\vert \xi _{k}\right\vert ^{m}\right)
\left\Vert B\left( \lambda ,t\right) u\right\Vert _{E}.
\end{equation*}

In view of estimate $\left( 2.7\right) $ and by Condition 2.1 we get from
the above inequality 
\begin{equation*}
\left\Vert \sigma _{\alpha }\left( t,\lambda ,\xi \right) u\right\Vert
_{E}\leq C_{\alpha }\left\Vert u\right\Vert _{E}.
\end{equation*}

\bigskip\ So, we obtain that the operator functions $\sigma \left( t,\lambda
,\xi \right) $ and $\sigma _{\alpha }\left( t,\lambda ,\xi \right) $ are
uniformly bounded, i.e., 
\begin{equation}
\left\Vert \sigma \left( t,\lambda ,\xi \right) \right\Vert _{B\left(
E\right) }\leq C,\left\Vert \sigma _{\alpha }\left( t,\lambda ,\xi \right)
\right\Vert _{B\left( E\right) }\leq C_{\alpha }.  \tag{2.8}
\end{equation}

Due to $R$-positivity of $A,$ by $\left( 2.8\right) $ and by Kahane's
contraction principle $\left[ \text{6, Lemma 3.5}\right] $\ we obtain that
the set 
\begin{equation*}
\left\{ \sigma \left( t,\lambda ,\xi \right) ;\text{ }\xi \in
R^{n}\backslash \left\{ 0\right\} \right\}
\end{equation*}%
is uniformly $R$-bounded, i.e., 
\begin{equation*}
\sup\limits_{t,\lambda }R\left\{ \sigma \left( t,\lambda ,\xi \right) ;\text{
}\xi \in R^{n}\backslash \left\{ 0\right\} \right\} \leq M_{0}.
\end{equation*}

In a similar way we obtain%
\begin{equation}
R\left( \left\{ \left\vert \xi \right\vert ^{\left\vert \beta \right\vert
}D_{\xi }^{\beta }\sigma \left( t,\lambda ,\xi \right) :\xi \in
R^{n}\backslash \left\{ 0\right\} \right\} \right) \leq M  \tag{2.9}
\end{equation}

for 
\begin{equation*}
\beta =\left( \beta _{1},\beta _{2},...,\beta _{n}\right) ,\text{ }\beta
_{i}\in \left\{ 0,1\right\} .
\end{equation*}

Consider the following sets%
\begin{equation*}
\sigma ^{\beta }\left( t,\lambda ,\xi \right) =\left\{ \left\vert \xi
\right\vert ^{\left\vert \beta \right\vert }D_{\xi }^{\beta }\sigma \left(
t,\lambda ,\xi \right) :\xi \in R^{n}\backslash \left\{ 0\right\} \right\} ,%
\text{ }
\end{equation*}

\begin{equation*}
\sigma _{\alpha }^{\beta }\left( t,\lambda ,\xi \right) =\left\{ \left\vert
\xi \right\vert ^{\left\vert \beta \right\vert }D_{\xi }^{\beta }\sigma
_{\alpha }\left( t,\lambda ,\xi \right) :\xi \in R^{n}\backslash \left\{
0\right\} \right\} ,\text{ }
\end{equation*}

\begin{equation*}
\beta =\left( \beta _{1},\beta _{2},...,\beta _{n}\right) ,\text{ }\beta
_{i}\in \left\{ 0,1\right\} .
\end{equation*}%
In view of the $R-$positivity properties of operator $A$ and due to Kahane's
contraction, addition and product properties of the collection of $R$%
-bounded operators (see e.g. $\left[ \text{7, 23}\right] $) and by $\left(
2.9\right) $ for all $\left\{ \xi ^{\left( j\right) }\right\} \in R^{n}$, $%
\left\{ \sigma _{\alpha }^{\beta }\left( t,\lambda ,\xi ^{\left( j\right)
}\right) \right\} $, $j=1,2,...,\mu ,$ $u_{1,}u_{2},...,u_{\mu }\in E$ and
independent symmetric $\left\{ -1,1\right\} -$valued random variables $%
r_{j}\left( y\right) $, $\mu \in \mathbb{N}$ we obtain the following uniform
estimate

\begin{equation*}
\int\limits_{\Omega }\left\Vert \sum\limits_{j=1}^{\mu }r_{j}\left( y\right)
\sigma _{\alpha }^{\beta }\left( t,\lambda ,\xi ^{\left( j\right) }\right)
u_{j}\right\Vert _{E}dy\leq
\end{equation*}

\begin{equation*}
C\int\limits_{\Omega }\left\Vert \sum\limits_{j=1}^{\mu }\sigma ^{\beta
}\left( t,\lambda ,\xi ^{\left( j\right) }\right) r_{j}\left( y\right)
u_{j}\right\Vert _{E}dy\leq C\int\limits_{\Omega }\left\Vert
\sum\limits_{j=1}^{\mu }r_{j}\left( y\right) u_{j}\right\Vert _{E}dy,
\end{equation*}%
i.e.,

\begin{equation*}
R\left( \left\{ \xi ^{\beta }D_{\xi }^{\beta }\sigma _{\alpha }\left(
t,\lambda ,\xi \right) :\xi \in R^{n}\backslash \left\{ 0\right\} \right\}
\right) \leq M_{\beta }.
\end{equation*}%
Hence, we infer that the operator-valued functions $\sigma \left( t,\lambda
,\xi \right) $ and $\sigma _{\alpha }\left( t,\lambda ,\xi \right) $\ are
uniform $R-$bounded multipliers and it's $R-$bounds are independent of $t$
and $\lambda $. By virtue of\ Preposition A$_{0}$, the operator-valued
functions $\sigma \left( t,\lambda ,\xi \right) \ $and $\sigma _{\alpha
}\left( t,\lambda ,\xi \right) $ are uniform collections of Fourier
multipliers in $L_{p}\left( R^{n};E\right) .$ So, we obtain that for all \ $%
f\in L_{p}\left( R^{n};E\right) $ there is a unique solution of equation $%
\left( 2.1\right) $ and estimate $\left( 2.3\right) $ holds$.$

Let $O_{t}$ denote the operator in $X=L_{p}\left( R^{n};E\right) $ generated
by problem $\left( 2.1\right) $ for $\lambda =0$, i.e., 
\begin{equation*}
D\left( O_{t}\right) \subset H_{p}^{m}\left( R^{n};E\left( A\right)
,E\right) ,\text{ }O_{t}u=P_{t}\left( D\right) u+Au.
\end{equation*}

Theorem 2.1 and the definition of the space $H_{p}^{m}\left( R^{n};E\left(
A\right) ,E\right) $ imply the following result:

\textbf{Result 2.1. }Assume all conditions of Theorem 2.1 are satisfied.
Then there are positive constants $C_{1}$ and $C_{2}$ so that 
\begin{equation*}
C_{1}\left\Vert O_{t}u\right\Vert _{X}\leq \left\Vert u\right\Vert
_{H_{p,t}^{m}\left( R^{n};E\left( A\right) ,E\right) }\leq C_{2}\left\Vert
O_{t}u\right\Vert _{X}
\end{equation*}

for $u\in Y.$ Indeed, if we put $\lambda =1$ in $\left( 2.3\right) ,$ by
Theorem 2.1 we get 
\begin{equation}
\dsum\limits_{\left\vert \alpha \right\vert \leq m}t\left( \alpha \right)
\left\Vert D^{\alpha }u\right\Vert _{X}+\left\Vert Au\right\Vert _{X}\leq
C\left\Vert O_{t}u\right\Vert _{X}  \tag{2.10}
\end{equation}

for $u\in Y$. Due to the closedness of $A$ and by the differential
properties of the Fourier transform we have{\large 
\begin{equation*}
\left\Vert Au\right\Vert _{X}=\left\Vert F^{-1}A\hat{u}\right\Vert _{X},%
\text{ }\left\Vert D^{\alpha }u\right\Vert _{X}=\left\Vert F^{-1}\xi
^{\alpha }\hat{u}\right\Vert _{X}.
\end{equation*}%
}

So, in view of estimate $\left( 2.10\right) $ and by definition of $Y$ we
obtain%
\begin{equation*}
\left\Vert u\right\Vert _{H_{p,t}^{m}\left( R^{n};E\left( A\right) ,E\right)
}\leq C_{2}\left\Vert O_{t}u\right\Vert _{X}.
\end{equation*}

\bigskip\ The first inequality is equivalent to the following estimate%
\begin{equation*}
\left\Vert F^{-1}A\hat{u}\right\Vert _{X}+\left\Vert F^{-1}P_{t}\left( \xi
\right) \hat{u}\right\Vert _{X}\leq
\end{equation*}

\begin{equation*}
C\left\{ \left\Vert F^{-1}A\hat{u}\right\Vert _{X}+\left\Vert F^{-1}\left[
1+\left( \dsum\limits_{k=1}^{n}t_{k}^{\frac{2}{m}}\xi _{k}^{2}\right) ^{%
\frac{1}{2}}\right] ^{m}\hat{u}\right\Vert _{X}\right\} .
\end{equation*}%
So, it suffices to show that the operator functions%
\begin{equation*}
A\left\{ A+\left[ 1+\left( \dsum\limits_{k=1}^{n}t_{k}^{\frac{2}{m}}\xi
_{k}^{2}\right) ^{\frac{1}{2}}\right] ^{m}\right\} ^{-1},\text{ }t\left(
\alpha \right) \xi ^{\alpha }\left[ 1+\left( \dsum\limits_{k=1}^{n}t_{k}^{%
\frac{2}{m}}\xi _{k}^{2}\right) ^{\frac{1}{2}}\right] ^{-m}
\end{equation*}%
are uniform Fourier multipliers in $X.$ This fact is proved in a similar way
as in the proof of Theorem 2.1.

From Theorem 2.1 we have:

\textbf{Result 2.2. }Assume all conditions of Theorem 2.1 hold. Then, for
all $\lambda \in S_{\varphi }$ the resolvent of operator $O_{t}$ exists and
the following sharp uniform estimate holds 
\begin{equation}
\sum\limits_{\left\vert \alpha \right\vert \leq m}t\left( \alpha \right)
\left\Vert D^{\alpha }\left( O_{t}+\lambda \right) ^{-1}\right\Vert
_{B\left( X\right) }+\left\Vert A\left( O_{t}+\lambda \right)
^{-1}\right\Vert _{B\left( X\right) }\leq C.  \tag{2.11}
\end{equation}

Indeed, we infer from Theorem 2.1 that the operator $O_{t}+\lambda $ has a
bounded inverse from $X$ to $Y.$ So, the solution\ $u$\ of equation $\left(
2.1\right) $ can be expressed as $u\left( x\right) =\left( O_{t}+\lambda
\right) ^{-1}f$ for all $f\in X.$ Then estimate $\left( 2.3\right) $ implies
the estimate $\left( 2.11\right) .$

\textbf{Result 2.3. }Theorem 2.1 particularly implies that the operator $%
O_{t}$ is positive in $X.$ Then the operators $O_{t}^{\sigma }$ are
generators of analytic semigroups in $X$ for $\sigma \leq \frac{1}{2}$ (see
e.g. $\left[ \text{20, \S 1.14.5}\right] $)$.$

\bigskip Now consider the problem $\left( 1.1\right) .$ By using Theorem 2.1
and the perturbation theory of linear operators we have the following

\textbf{Theorem 2.2.} Assume all conditions of Theorem 2.1 are satisfied.
Suppose $A_{\alpha }\left( x\right) A^{-\left( 1-\frac{\left\vert \alpha
\right\vert }{m}-\mu \right) }\in L_{\infty }\left( R^{n};B\left( E\right)
\right) $ for $\mu \in \left( 0,1-\frac{\left\vert \alpha \right\vert }{m}%
\right) $. Then for $f\in L_{p}\left( R^{n};E\right) $, $\lambda \in
S_{\varphi _{2}},$ $0\leq \varphi _{2}<\pi -\varphi _{1}$, $\varphi
_{1}+\varphi _{2}\leq \varphi $ and for sufficiently large $\left\vert
\lambda \right\vert $\ there is a unique solution $u$ of the equation $%
\left( 1.1\right) $ belonging to $Y$ and the following coercive uniform
estimate holds 
\begin{equation}
\dsum\limits_{\left\vert \alpha \right\vert \leq m}t\left( \alpha \right)
\left\vert \lambda \right\vert ^{1-\frac{\left\vert \alpha \right\vert }{m}%
}\left\Vert D^{\alpha }u\right\Vert _{L_{p}\left( R^{n};E\right)
}+\left\Vert Au\right\Vert _{L_{p}\left( R^{n};E\right) }\leq C\left\Vert
f\right\Vert _{L_{p}\left( R^{n};E\right) }.  \tag{2.12}
\end{equation}

\textbf{Proof. }It is clear that $Q_{t}=O_{t}+L_{t},$ where $O_{t}$ is the
operator in $L_{p}\left( R^{n};E\right) $\ generated by problem $\left(
2.1\right) $ for $\lambda =0$ and 
\begin{equation*}
\ L_{t}u=\sum\limits_{\left\vert \alpha \right\vert <m}t\left( \alpha
\right) A_{\alpha }\left( x\right) D^{\alpha }u,\text{ }u\in Y.
\end{equation*}

In view of the condition on $A_{\alpha }\left( x\right) $ and by the Theorem
A$_{1}$ for$\ u\in Y$\ we have\ \ \ 

\begin{equation}
\left\Vert L_{t}u\right\Vert _{X}\leq \sum\limits_{\left\vert \alpha
\right\vert <m}t\left( \alpha \right) \left\Vert A_{\alpha }\left( x\right)
D^{\alpha }u\right\Vert _{X}\leq  \notag
\end{equation}%
\begin{equation*}
C_{\mu }\sum\limits_{\left\vert \alpha \right\vert <m}t\left( \alpha \right)
\left\Vert A^{1-\frac{\left\vert \alpha \right\vert }{m}-\mu }D^{\alpha
}u\right\Vert _{X}\leq
\end{equation*}%
\begin{equation}
C_{\mu }\left[ h^{\mu }\left\Vert u\right\Vert _{H_{p,t}^{m}\left(
R^{n};E\left( A\right) ,E\right) }+h^{-\left( 1-\mu \right) }\left\Vert
u\right\Vert _{X}\right] .  \tag{2.13}
\end{equation}%
Then from estimates $\left( 2.12\right) ,$ $\left( 2.13\right) $ and for\ $%
u\in Y$ we\ obtain

\begin{equation}
\left\Vert L_{t}u\right\Vert _{X}\leq C\left[ h^{\mu }\left\Vert
O_{t}u\right\Vert _{X}+h^{-\left( 1-\mu \right) }\left\Vert u\right\Vert
_{X}.\right]  \tag{2.14}
\end{equation}%
Since \ $\left\Vert u\right\Vert _{X}=\frac{1}{\lambda }\left\Vert \left(
O_{t}+\lambda \right) u+L_{t}u\right\Vert _{X}$ \ for $\lambda \in
S_{\varphi _{2}}.$ Hence, for $u\in Y$ we get 
\begin{equation}
\left\Vert u\right\Vert _{X}\leq \frac{1}{\left\vert \lambda \right\vert }%
\left[ \left\Vert \left( O_{t}+\lambda \right) u\right\Vert _{X}+\left\Vert
O_{t}u\right\Vert _{X}\right] \leq  \tag{2.15}
\end{equation}%
\begin{equation*}
\leq \frac{1}{\left\vert \lambda \right\vert }\left\Vert \left(
O_{t}+\lambda \right) u\right\Vert _{X}+\frac{1}{\left\vert \lambda
\right\vert }\left[ \sum\limits_{\left\vert \alpha \right\vert <m}t\left(
\alpha \right) \left\Vert D^{\alpha }u\right\Vert _{X}+\left\Vert
Au\right\Vert _{X}\right] .
\end{equation*}%
From estimates $\left( 2.13\right) -\left( 2.15\right) $ for\ $u\in Y$ we\
obtain 
\begin{equation}
\left\Vert L_{t}u\right\Vert _{X}\leq Ch^{\mu }\left\Vert \left(
O_{t}+\lambda \right) u\right\Vert _{X}+C_{1}\left\vert \lambda \right\vert
^{-1}h^{-\left( 1-\mu \right) }\left\Vert \left( O_{t}+\lambda \right)
u\right\Vert _{X}.  \tag{2.16}
\end{equation}%
Then choosing $h$ and $\lambda $ such that $Ch^{\mu }<1,$ \ $C_{1}\left\vert
\lambda \right\vert ^{-1}h^{-\left( 1-\mu \right) }<1$ from $\left(
2.16\right) $ we obtain that 
\begin{equation}
\ \ \left\Vert L_{t}\left( O_{t}+\lambda \right) ^{-1}\right\Vert _{B\left(
X\right) }<1.  \tag{2.17}
\end{equation}%
From Theorem 2.1 and $\left( 2.17\right) $ we get that the operator $\left(
Q_{t}+\lambda \right) $ has a bounded inverse in $X.$ Moreover, it is clear
that 
\begin{equation}
\left( Q_{t}+\lambda \right) ^{-1}=\left[ I+L_{t}\left( O_{t}+\lambda
\right) ^{-1}\right] \left( O_{t}+\lambda \right) ,  \tag{2.18}
\end{equation}%
\ where $I$ is an identity operator in $X.$ Using relation $\left(
2.18\right) $, estimates $\left( 2.3\right) $, $\left( 2.17\right) $ and
perturbation theory of linear operators, we obtain that the operator $%
Q_{t}+\lambda $\ has a bounded inverse from $X$ into $Y$ and the estimate $%
\left( 2.12\right) $ holds.

\begin{center}
\textbf{3. The Cauchy problem for parabolic PsDOE with parameter}
\end{center}

In this section, we shall consider the following Cauchy problem for the
parabolic PsDO equation 
\begin{equation}
\frac{\partial u}{\partial y}+P_{t}\left( D\right) u+Au=f\left( y,x\right) ,%
\text{ }u(0,x)=0,  \tag{3.1}
\end{equation}%
where $P_{t}\left( D\right) $ is the pseudo-differential operator defined by 
$\left( 2.2\right) $ and $A$ is a linear operator in $E$, $t=\left(
t_{1},t_{2},...,t_{n}\right) $, $t_{k}$ are positive parameters.

In this section, by applying Theorem 2.1 we establish the maximal regularity
of the problem $\left( 3.1\right) $ in $E-$valued mixed $L_{\mathbf{p}}$
spaces, where $\mathbf{p=}\left( p_{1},p\right) .$

Let $O_{t}$ denote the operator generated by problem $\left( 2.1\right) .$
For this aim we need the following result:

\textbf{Theorem 3.1. }Suppose Condition 2. 1 hold, $E$ is an UMD space and
the operator $A$ is $R$-positive in $E$ with respect to $\varphi $ with $%
0\leq \varphi <\pi -\varphi _{1}$. Then operator $O_{t}$ is uniformly $R$%
-positive in $L_{p}\left( R^{n};E\right) .$

\bigskip \textbf{Proof. }From Result 2.3 we obtain that the operator $O_{t}$
is positive in $X=L_{p}\left( R^{n};E\right) $. We have to prove the $R$%
-boundedness of the set 
\begin{equation*}
\sigma \left( t,\lambda ,\xi \right) =\left\{ \lambda \left( O_{t}+\lambda
\right) ^{-1}:\lambda \in S_{\varphi }\right\} .
\end{equation*}%
From the proof of Theorem 2.1 we have 
\begin{equation*}
\lambda \left( O_{t}+\lambda \right) ^{-1}f=F^{-1}\Phi \left( t,\xi ,\lambda
\right) \hat{f}\text{, }f\in L_{p}\left( R^{n};E\right) ,
\end{equation*}%
where 
\begin{equation*}
\Phi \left( t,\xi ,\lambda \right) =\lambda \left( A+P_{t}\left( \xi \right)
+\lambda \right) ^{-1}.
\end{equation*}%
By reasoning as in the proof of Theorem 2.1, we obtain the following uniform
estimate%
\begin{equation*}
\left\Vert \Phi \left( t,\xi ,\lambda \right) \right\Vert _{B\left( E\right)
}\leq \left\vert \lambda \right\vert \left\Vert \left( A+P_{t}\left( \xi
\right) +\lambda \right) ^{-1}\right\Vert _{B\left( E\right) }\leq C.
\end{equation*}%
By definition of $R$-boundedness, it suffices to show that the operator
function $\Phi \left( t,\xi ,\lambda \right) $ ( which depends on variable $%
\lambda $ and parameters $\xi ,$ $t$ ) is a multiplier in $L_{p}\left(
R^{n};E\right) $ uniformly\ with respect to $\zeta $ and $t.$ Indeed, by
reasoning as in Theorem 2.1 we can easily show that $\Phi \left( t,\xi
,\lambda \right) $ is a uniform multiplier in $L_{p}\left( R;E\right) .$
Then,\ by the definition of a $R$-bounded set we have 
\begin{equation*}
\int\limits_{0}^{1}\left\Vert \sum\limits_{j=1}^{m}r_{j}\left( y\right)
\lambda _{j}\left( O+\lambda _{j}\right) ^{-1}f_{j}\right\Vert
_{X}dy=\int\limits_{0}^{1}\left\Vert \sum\limits_{j=1}^{m}r_{j}\left(
y\right) F^{-1}\Phi \left( t,\xi ,\lambda _{j}\right) \hat{f}_{j}\right\Vert
_{X}dy
\end{equation*}

\begin{equation*}
=\int\limits_{0}^{1}\left\Vert F^{-1}\sum\limits_{j=1}^{m}r_{j}\left(
y\right) \Phi \left( t,\xi ,\lambda _{j}\right) \hat{f}_{j}\right\Vert
_{X}dy\leq C\int\limits_{0}^{1}\left\Vert \sum\limits_{j=1}^{m}r_{j}\left(
y\right) f_{j}\right\Vert _{X}dy
\end{equation*}%
for all $\xi \in R$, $\lambda _{1},\lambda _{2}...\lambda _{m}\in S_{\varphi
},$ $f_{1,}f_{2},...,f_{m}\in X$, $m\in N$, where $\left\{ r_{j}\right\} $
is a sequence of independent symmetric $\left\{ -1,1\right\} -$valued random
variables on $\left[ 0,1\right] $. Hence, the set $\Phi \left( t,\xi
,\lambda \right) $ is uniformly $R$-bounded.

Let $E$ be a Banach space. For $\mathbf{p=}\left( p\text{, }p_{1}\right) ,$ $%
Z=L_{\mathbf{p}}\left( R_{+}^{n+1};E\right) $ will denote the space of all $%
\mathbf{p}$-summable $E$-valued\ functions with mixed norm (see e.g. $\left[ 
\text{2, \S\ 4}\right] $ for the complex-valued case), i.e., the space of
all measurable $E$-valued functions $f$ defined on $R_{+}^{n+1}$, for which 
\begin{equation*}
\left\Vert f\right\Vert _{L_{\mathbf{p}}\left( R_{+}^{n+1};E\right) }=\left(
\int\limits_{R^{n}}\left( \int\limits_{R_{+}}\left\Vert f\left( y,x\right)
\right\Vert _{E}^{p}dx\right) ^{\frac{p_{1}}{p}}dy\right) ^{\frac{1}{p_{1}}%
}<\infty .
\end{equation*}

Let $E$ be a Banach space and $A$ be a positive operator in $E.$ Suppose, $l$
is a positive integer number. $W_{p}^{l}\left( R_{+};E\left( A\right)
,E\right) $ denotes the space of all functions $u\in L_{p}\left(
R_{+};E\left( A\right) \right) $ possessing the generalized derivatives $%
u^{\left( l\right) }\in L_{p}\left( R_{+};E\right) $ with the norm 
\begin{equation*}
\ \left\Vert u\right\Vert _{W_{p}^{l}\left( R_{+};E\left( A\right) ,E\right)
}=\left\Vert Au\right\Vert _{L_{p}\left( R_{+};E\right) }+\left\Vert
u^{\left( l\right) }\right\Vert _{L_{p}\left( R_{+};E\right) }.
\end{equation*}%
Let $m$ be a positive number. $W_{\mathbf{p}}^{1,m}\left(
R_{+}^{n+1};E\left( A\right) ,E\right) $ denotes the space of all functions $%
u\in L_{\mathbf{p}}\left( R_{+}^{n+1};E\left( A\right) \right) $ possessing
the generalized derivative $D_{y}u=\frac{\partial u}{\partial y}\in Z$ with
respect to $y$\ and fractional derivatives $D_{x}^{\alpha }u\in Z$ with
respect to $x$ for $\left\vert \alpha \right\vert \leq m$ with the norm 
\begin{equation*}
\ \left\Vert u\right\Vert _{W_{\mathbf{p}}^{1,m}\left( R_{+}^{n+1};E\left(
A\right) ,E\right) }=\left\Vert Au\right\Vert _{Z}+\left\Vert \frac{du}{dy}%
\right\Vert _{Z}+\dsum\limits_{\left\vert \alpha \right\vert \leq
m}\left\Vert D_{x}^{\alpha }u\right\Vert _{Z},
\end{equation*}

where $u=u\left( y,x\right) .$

Now, we are ready to state the main result of this section.

\textbf{Theorem 3.2.}\ Assume the conditions of Theorem 2.1 hold for $%
\varphi \in \left( \frac{\pi }{2},\pi \right) $. Then for $f\in Z$ problem $%
\left( 3.1\right) $ has a unique solution%
\begin{equation*}
u\in W_{\mathbf{p}}^{1,m}\left( R_{+}^{n+1};E\left( A\right) ,E\right)
\end{equation*}%
satisfying the following unform coercive estimate 
\begin{equation*}
\left\Vert \frac{du}{dy}\right\Vert _{Z}+\left\Vert P_{t}\left( D\right)
u\right\Vert _{Z}+\left\Vert Au\right\Vert _{Z}\leq C\left\Vert f\right\Vert
_{Z}.
\end{equation*}

\textbf{Proof.} By definition of $X=L_{p}\left( R^{n};E\right) $ and mixed
space $L_{\mathbf{p}}\left( R_{+}^{n+1};E\right) ,$ $\mathbf{p=}\left( p%
\text{, }p_{1}\right) $, we have

\begin{equation*}
\left\Vert u\right\Vert _{L_{p_{1}}\left( R_{+};X\right) }=\left(
\dint\limits_{R_{+}}\left\Vert u\left( y\right) \right\Vert
_{X}^{p_{1}}dy\right) ^{\frac{1}{p_{1}}}=\left(
\dint\limits_{R_{+}}\left\Vert u\left( y\right) \right\Vert _{L_{p}\left(
R^{n};E\right) }^{p_{1}}dy\right) ^{\frac{1}{p_{1}}}=
\end{equation*}%
\begin{equation*}
\left( \int\limits_{R^{n}}\left( \int\limits_{R_{+}}\left\Vert u\left(
y,x\right) \right\Vert _{E}^{p}dy\right) ^{\frac{p_{1}}{p}}dx\right) ^{\frac{%
1}{p_{1}}}=\left\Vert u\right\Vert _{Z}.
\end{equation*}%
Moreover, by definition of the space $W_{p}^{m}\left( R_{+};E\left( A\right)
,E\right) $ and by Result 2.1 we obtain 
\begin{equation}
\left\Vert u\right\Vert _{W_{p_{1}}^{1}\left( R_{+};D\left( O_{t}\right)
,X\right) }=\left\Vert O_{t}u\right\Vert _{L_{p}\left( R_{+};X\right)
}+\left\Vert u^{\prime }\right\Vert _{L_{p}\left( R_{+};X\right) }= 
\tag{3.2}
\end{equation}%
\begin{equation*}
\left\Vert Au\right\Vert _{Z}+\left\Vert D_{y}u\right\Vert
_{Z}+\dsum\limits_{\left\vert \alpha \right\vert \leq m}\left\Vert
D_{x}^{\alpha }u\right\Vert _{Z}=\left\Vert u\right\Vert _{W_{\mathbf{p}%
}^{1,m}\left( R_{+}^{n+1};E\left( A\right) ,E\right) .}
\end{equation*}

Hence, we deduced from the above equalities that, 
\begin{equation*}
Z=L_{\mathbf{p}}\left( R_{+}^{n+1};E\right) =L_{p_{1}}\left( R_{+};X\right) 
\text{, }W_{\mathbf{p}}^{1,m}\left( R_{+}^{n+1};E\left( A\right) ,E\right) =
\end{equation*}%
\begin{equation*}
W_{p_{1}}^{1}\left( R_{+};D\left( O_{t}\right) ,X\right) .
\end{equation*}%
Therefore, the problem $\left( 3.1\right) $\ can be expressed as the
following Cauchy problem for the abstract parabolic equation 
\begin{equation}
\frac{du}{dy}+O_{t}u\left( y\right) =f\left( y\right) ,\text{ }u\left(
0\right) =0,\text{ }y\in R_{+}.  \tag{3.3}
\end{equation}

By virtue of $\left[ \text{1, Theorem 4.5.2}\right] $, the condition $E\in
UMD$ implies $X\in UMD$ for $p\in \left( 1,\infty \right) $. Then due to the 
$R-$positivity of $O_{t},$ by virtue of [23, Theorem 4.2] we obtain that for 
$f\in L_{p_{1}}\left( R_{+};X\right) $ equation $\left( 3.3\right) $ has a
unique solution $u\in W_{p_{1}}^{1}\left( R_{+};D\left( O_{t}\right)
,X\right) $ satisfying the following estimate 
\begin{equation*}
\left\Vert \frac{du}{dy}\right\Vert _{L_{p_{1}}\left( R_{+};X\right)
}+\left\Vert O_{t}u\right\Vert _{L_{p_{1}}\left( R_{+};X\right) }\leq
C\left\Vert f\right\Vert _{L_{p_{1}}\left( R_{+};X\right) }.
\end{equation*}

From the Theorem 2.1, relation $\left( 3.2\right) $ and from the above
estimate we get the assertion$.$

\begin{center}
\textbf{4. BVP for Anisotropic PsDE }
\end{center}

In this section, the maximal regularity properties of the anisotropic PsDE
are studied.

Let $\tilde{\Omega}=\Omega \times R^{n}$, where $\Omega \subset R^{\mu }$ is
an open connected set with compact $C^{2l}-$boundary $\partial \Omega $.
Consider the BVP for the pseudo-differential equation

\begin{equation}
P_{t}\left( D\right) u+\sum\limits_{\left\vert \alpha \right\vert \leq
2l}\left( b_{\alpha }D_{y}^{\alpha }+\lambda \right) u=f\left( x,y\right) ,%
\text{ }y\in \Omega ,  \tag{4.1}
\end{equation}%
\ \ \ 

\begin{equation}
B_{j}u=\sum\limits_{\left\vert \beta \right\vert \leq l_{j}}\ b_{j\beta
}\left( y\right) D_{y}^{\beta }u\left( x,y\right) =0\text{, }x\in R^{n}\text{%
,}  \tag{4.2}
\end{equation}

\begin{equation*}
\text{ }y\in \partial \Omega ,\text{ }j=1,2,...,l,
\end{equation*}%
where $u=u\left( x,y\right) ,\ P_{t}\left( D\right) $ is the pseudo
differential operator defined by $\left( 2.1\right) $ with respect to $x$
and 
\begin{equation*}
D_{j}=-i\frac{\partial }{\partial y_{j}}\text{, }y=\left( y_{1},...,y_{\mu
}\right) \text{, }b_{\alpha }=b_{\alpha }\left( y\right) \text{,}
\end{equation*}

where $\alpha =\left( \alpha _{1},\alpha _{2},...,\alpha _{\mu }\right) $, $%
\beta =\left( \beta _{1},\beta _{2},...,\beta _{\mu }\right) $ are
nonnegative integer numbers, $t=\left( t_{1},t_{2},...,t_{n}\right) $ and $%
t_{k}$ are positive parameters.

If \ $\tilde{\Omega}=R^{n}\times \Omega $, $\mathbf{p=}\left( p_{1},p\right) 
$, $L_{\mathbf{p}}\left( \tilde{\Omega}\right) $ will denote the space of
all $\mathbf{p}$-summable scalar-valued\ functions with mixed norm ( see e.g.%
$\left[ \text{2, \S\ 4}\right] $ ), i.e., the space of all measurable
functions $f$ defined on $\tilde{\Omega}$, for which 
\begin{equation*}
\left\Vert f\right\Vert _{L_{\mathbf{p}}\left( \tilde{\Omega}\right)
}=\left( \int\limits_{R^{n}}\left( \int\limits_{\Omega }\left\vert f\left(
x,y\right) \right\vert ^{p_{1}}dx\right) ^{\frac{p}{p_{1}}}dy\right) ^{\frac{%
1}{p}}<\infty .
\end{equation*}%
Analogously, $W_{\mathbf{p}}^{m,2l}\left( \tilde{\Omega}\right) $ denotes
the anisotropic fractional Sobolev space with corresponding mixed norm,
i.e., $W_{\mathbf{p}}^{m,2l}\left( \tilde{\Omega}\right) $ denotes the space
of all functions $u\in L_{\mathbf{p}}\left( \tilde{\Omega}\right) $
possessing the fractional derivatives $D_{x}^{\alpha }u\in L_{\mathbf{p}%
}\left( \tilde{\Omega}\right) $ with respect to $x$ for $\left\vert \alpha
\right\vert \leq m$ and generalized derivative $\frac{\partial ^{2l}u}{%
\partial y_{k}^{2l}}\in L_{\mathbf{p}}\left( \tilde{\Omega}\right) $ with
respect to $y$ with the norm 
\begin{equation*}
\ \left\Vert u\right\Vert _{W_{\mathbf{p}}^{m,2l}\left( \tilde{\Omega}%
\right) }=\dsum\limits_{\left\vert \alpha \right\vert \leq m}\left\Vert
D_{x}^{\alpha }u\right\Vert _{L_{\mathbf{p}}\left( \tilde{\Omega}\right)
}+\dsum\limits_{k=1}^{\mu }\left\Vert \frac{\partial ^{2l}u}{\partial
y_{k}^{2l}}\right\Vert _{L_{\mathbf{p}}\left( \tilde{\Omega}\right) }.
\end{equation*}

\ Let $Q$ denote the operator generated by problem $\left( 4.1\right)
-\left( 4.2\right) $, i.e., 
\begin{equation*}
D\left( Q\right) =W_{\mathbf{p}}^{m,2l}\left( \tilde{\Omega},B_{j}\right)
=\left\{ u:u\in W_{\mathbf{p}}^{m,2l}\left( \tilde{\Omega}\right) ,\text{ }%
B_{j}u=0,\text{ }j=1,2,...l\right\} ,
\end{equation*}%
\begin{equation*}
Qu=P_{t}\left( D\right) u+\sum\limits_{\left\vert \alpha \right\vert \leq
2l}b_{\alpha }D_{y}^{\alpha }u.
\end{equation*}

Let $\xi ^{\prime }=\left( \xi _{1},\xi _{2},...,\xi _{\mu -1}\right) \in
R^{\mu -1},$ $\alpha ^{\prime }=\left( \alpha _{1},\alpha _{2},...,\alpha
_{\mu -1}\right) \in Z^{\mu }$ and 
\begin{equation*}
\text{ }A\left( y_{0},\xi ^{\prime },D_{y}\right) =\sum\limits_{\left\vert
\alpha ^{\prime }\right\vert +j\leq 2l}a_{\alpha ^{\prime }}\left(
y_{0}\right) \xi _{1}^{\alpha _{1}}\xi _{2}^{\alpha _{2}}...\xi _{\mu
-1}^{\alpha _{\mu -1}}D_{\mu }^{j}\text{ for }y_{0}\in \bar{G}
\end{equation*}%
\begin{equation*}
B_{j}\left( y_{0},\xi ^{\prime },D_{y}\right) =\sum\limits_{\left\vert \beta
^{\prime }\right\vert +j\leq l_{j}}b_{j\beta ^{\prime }}\left( y_{0}\right)
\xi _{1}^{\beta _{1}}\xi _{2}^{\beta _{2}}...\xi _{\mu -1}^{\beta _{\mu
-1}}D_{\mu }^{j}\text{ for }y_{0}\in \partial G.
\end{equation*}

\textbf{Condition 4.1. }Let the following conditions be satisfied;

(1) $b_{\alpha }\in C\left( \bar{\Omega}\right) $ for each $\left\vert
\alpha \right\vert =2l$ and $b_{\alpha }\in L_{\infty }\left( \Omega \right)
+L_{r_{k}}\left( \Omega \right) $ for each $\left\vert \alpha \right\vert
=k<2l$ with $r_{k}\geq p_{1}$, $p_{1}\in \left( 1,\infty \right) $ and $2l-k>%
\frac{l}{r_{k}};$

(2) $b_{j\beta }\in C^{2l-l_{j}}\left( \partial \Omega \right) $ for each $j$%
, $\beta $, $l_{j}<2l$, $p\in \left( 1,\infty \right) ,$ $\lambda \in
S_{\varphi },$ $\varphi \in \lbrack 0,\pi );$

(3) for $y\in \bar{\Omega}$, $\xi \in R^{\mu }$, $\sigma \in S_{\varphi
_{0}} $, $\varphi _{0}\in \left( 0,\frac{\pi }{2}\right) $, $\left\vert \xi
\right\vert +\left\vert \sigma \right\vert \neq 0$ let $\sigma
+\sum\limits_{\left\vert \alpha \right\vert =2l}b_{\alpha }\left( y\right)
\xi ^{\alpha }\neq 0;$

(4) for each $y_{0}\in \partial \Omega $ local BVP in local coordinates
corresponding to $y_{0}$ 
\begin{equation*}
\lambda +A\left( y_{0},\xi ^{\prime },D_{y}\right) \vartheta \left( y\right)
=0,
\end{equation*}

\begin{equation*}
B_{j}\left( y_{0},\xi ^{\prime },D_{y}\right) \vartheta \left( 0\right)
=h_{j}\text{, }j=1,2,...,l
\end{equation*}%
has a unique solution $\vartheta \in C_{0}\left( \mathbb{R}_{+}\right) $ for
all $h=\left( h_{1},h_{2},...,h_{l}\right) \in \mathbb{C}^{l}$ and for $\xi
^{^{\prime }}\in R^{n-1}$.

\textbf{\ }Suppose $\nu =\left( \nu _{1},\nu _{2},...,\nu _{n}\right) $ are
nonnegative real numbers. In this section, we present the following result:

\textbf{Theorem 4.1}. Assume\textbf{\ }Condition 2.1 and Condition 4.1 are
satisfied. Then for $\ f\in L_{\mathbf{p}}\left( \tilde{\Omega}\right) $, $%
\lambda \in S_{\varphi },$ $\varphi \in \left( 0\right. ,\left. \pi \right] $
problem $\left( 4.1\right) -\left( 4.2\right) $ has a unique solution $u\in
W_{p}^{m,2l}\left( \tilde{\Omega}\right) $ and the following coercive
uniform estimate holds 
\begin{equation*}
\dsum\limits_{\left\vert \nu \right\vert \leq
m}\dprod\limits_{k=1}^{n}t_{k}^{\frac{\nu _{k}}{m}}\left\vert \lambda
\right\vert ^{1-\frac{\left\vert \nu \right\vert }{m}}\left\Vert D_{x}^{\nu
}u\right\Vert _{L_{\mathbf{p}}\left( \tilde{\Omega}\right)
}+\sum\limits_{\left\vert \alpha \right\vert \leq 2l}\left\Vert
D_{y}^{\alpha }u\right\Vert _{L_{\mathbf{p}}\left( \tilde{\Omega}\right)
}\leq C\left\Vert f\right\Vert _{L_{\mathbf{p}}\left( \tilde{\Omega}\right)
}.
\end{equation*}

\ \textbf{Proof.} Let $E=L_{p_{1}}\left( \Omega \right) $. It is known $%
\left[ 4\right] $\ that $L_{p_{1}}\left( \Omega \right) $ is an $UMD$ space
for $p_{1}\in \left( 1,\infty \right) .$ Consider the operator $A$ defined
by 
\begin{equation*}
D\left( A\right) =W_{p_{1}}^{2l}\left( \Omega ;B_{j}u=0\right) ,\text{ }%
Au=\sum\limits_{\left\vert \alpha \right\vert \leq 2l}b_{\alpha }\left(
x\right) D^{\alpha }u\left( y\right) .
\end{equation*}

Therefore, the problem $\left( 4.1\right) -\left( 4.2\right) $ can be
rewritten in the form of $\left( 2.1\right) $, where $u\left( x\right)
=u\left( x,.\right) ,$ $f\left( x\right) =f\left( x,.\right) $\ are
functions with values in $E=L_{p_{1}}\left( \Omega \right) .$ From $\left[ 
\text{6, Theorem 8.2}\right] $ we get that the following problem

\begin{equation}
\eta u\left( y\right) +\sum\limits_{\left\vert \alpha \right\vert \leq
2l}b_{\alpha }\left( y\right) D^{\alpha }u\left( y\right) =f\left( y\right) 
\text{, }  \tag{4.3}
\end{equation}%
\begin{equation*}
B_{j}u=\sum\limits_{\left\vert \beta \right\vert \leq l_{j}}\ b_{j\beta
}\left( y\right) D^{\beta }u\left( y\right) =0\text{, }j=1,2,...,l
\end{equation*}%
has a unique solution for $f\in L_{p_{1}}\left( \Omega \right) $ and arg $%
\eta \in S\left( \varphi _{1}\right) ,$ $\left\vert \eta \right\vert
\rightarrow \infty .$ Moreover, the operator $A$ generated\ by $\left(
4.3\right) $ is $R$-positive in $L_{p_{1}}.$ Then from Theorem 2.1 we obtain
the assertion.

\ \ \ \ \ 

\begin{center}
\textbf{5. The system of PsDE of infinite order }
\end{center}

Consider the following system of PsDEs of infinite order

\begin{equation}
P_{t}\left( D\right) u_{i}+\sum\limits_{j=1}^{N}\left( a_{ij}+\lambda
\right) u_{j}\left( x\right) =f_{i}\left( x\right) \text{, }x\in R^{n},\text{
}  \tag{5.1}
\end{equation}

\ \ \ 
\begin{equation*}
i=1,2,...,N,\text{ }N\in \left[ 1,\infty \right] ,
\end{equation*}%
where $P_{t}\left( D\right) $ is the pseudo-differential operator defined by 
$\left( 2.2\right) ,$ $t=\left( t_{1},t_{2},...,t_{n}\right) $ and $t_{k}$
are positive parameters. Let $a_{ij}$ be real numbers and%
\begin{equation*}
\text{ }l_{q}\left( A\right) =\left\{ u\in l_{q},\left\Vert u\right\Vert
_{l_{q}\left( A\right) }=\left\Vert Au\right\Vert _{l_{q}}=\right.
\end{equation*}

\begin{equation*}
\left( \sum\limits_{i=1}^{N}\left\vert \left( Au\right) _{i}\right\vert
^{q}\right) ^{\frac{1}{q}}=\left. \left( \sum\limits_{i=1}^{N}\left\vert
\sum\limits_{j=1}^{N}a_{ij}u_{j}\right\vert ^{q}\right) ^{\frac{1}{q}%
}<\infty \right\} ,
\end{equation*}%
\begin{equation*}
\text{ }u=\left\{ u_{j}\right\} ,\text{ }Au=\left\{
\dsum\limits_{j=1}^{N}a_{ij}u_{j}\right\} ,\text{ }i,\text{ }j=1,2,...N.
\end{equation*}

\textbf{Condition 5.1. }Let 
\begin{equation*}
a_{ij}=a_{ji}\text{, }\sum\limits_{i,j=1}^{N}a_{ij}\xi _{i}\xi _{j}\geq
C_{0}\left\vert \xi \right\vert ^{2},\text{ for }\xi \neq 0.
\end{equation*}

Let%
\begin{equation*}
f\left( x\right) =\left\{ f_{i}\left( x\right) \right\} _{1}^{N}\text{, }%
u=\left\{ u_{i}\left( x\right) \right\} _{1}^{N}.
\end{equation*}

\textbf{Theorem 5.1. }Assume\textbf{\ }Condition 2.1 and Condition 5.1 are
satisfied.\textbf{\ }Then, for $f\left( x\right) \in L_{p}\left(
R^{n};l_{q}\right) ,$ $\left\vert \arg \lambda \right\vert \leq \varphi $, $%
\varphi \in \left( 0\right. ,\left. \pi \right] $ and for sufficiently large 
$\left\vert \lambda \right\vert ,$ problem $\left( 5.1\right) $ has a unique
solution $u$ that belongs to the space $H_{p}^{m}\left( R^{n},l_{q}\left(
A\right) ,l_{q}\right) $ and the following uniform coercive estimate holds 
\begin{equation*}
\sum\limits_{\left\vert \alpha \right\vert \leq m}t\left( \alpha \right)
\left\vert \lambda \right\vert ^{1-\frac{\left\vert \alpha \right\vert }{m}}%
\left[ \int\limits_{R^{n}}\left( \sum\limits_{j=1}^{N}\left\vert D^{\alpha
}u_{j}\left( x\right) \right\vert ^{q}\right) ^{\frac{p}{q}}dx\right] ^{%
\frac{1}{p}}+
\end{equation*}

\begin{equation*}
\left[ \int\limits_{R^{n}}\left( \sum\limits_{i=1}^{N}\left\vert
\sum\limits_{j=1}^{N}a_{ij}u_{j}\right\vert ^{q}\right) ^{\frac{p}{q}}dx%
\right] ^{\frac{1}{p}}\leq C\left[ \int\limits_{R^{n}}\left(
\sum\limits_{i=1}^{N}\left\vert f_{i}\left( x\right) \right\vert ^{q}\right)
^{\frac{p}{q}}dx\right] ^{\frac{1}{p}}.
\end{equation*}%
\ \textbf{Proof.} Let $E=l_{q},$ $A$ be a matrix such that $A=\left[ a_{ij}%
\right] ,$ $i$, $j=1,2,...N.$ It is easy to see that 
\begin{equation*}
B\left( \lambda \right) =\lambda \left( A+\lambda \right) ^{-1}=\frac{%
\lambda }{D\left( \lambda \right) }\left[ A_{ji}\left( \lambda \right) %
\right] \text{, }i\text{, }j=1,2,...N,
\end{equation*}%
where $D\left( \lambda \right) =\det \left( A-\lambda I\right) $, $%
A_{ji}\left( \lambda \right) $ are entries of the corresponding adjoint
matrix of $A-\lambda I.$ Since the matrix $A$ is symmetric and positive
definite, it generates a positive operator in $l_{q}$ for $q\in \left(
1,\infty \right) .$ For all $u_{1,}u_{2},...,u_{\mu }\in l_{q}$, $\lambda
_{1},\lambda _{2},...,\lambda _{\mu }\in \mathbb{C}$ and independent
symmetric $\left\{ -1,1\right\} $-valued random variables $r_{k}\left(
y\right) $, $k=1,2,...,\mu ,$ $\mu \in \mathbb{N}$\ we have 
\begin{equation*}
\int\limits_{\Omega }\left\Vert \sum\limits_{k=1}^{\mu }r_{k}\left( y\right)
B\left( \lambda _{k}\right) u_{k}\right\Vert _{l_{q}}^{q}dy\leq
\end{equation*}%
\begin{equation*}
C\left\{ \int\limits_{\Omega }\sum\limits_{j=1}^{N}\left\vert
\sum\limits_{k=1}^{\mu }\sum\limits_{j=1}^{N}\frac{\lambda _{k}}{D\left(
\lambda _{k}\right) }A_{ji}\left( \lambda _{k}\right) r_{k}\left( y\right)
u_{ki}\right\vert ^{q}dy\right. \leq
\end{equation*}%
\begin{equation}
\sup\limits_{k,i}\sum\limits_{j=1}^{N}\left\vert \frac{\lambda _{k}}{D\left(
\lambda _{k}\right) }A_{ji}\left( \lambda _{k}\right) \right\vert
^{q}\int\limits_{\Omega }\left\vert \sum\limits_{k=1}^{\mu }r_{k}\left(
y\right) u_{kj}\right\vert ^{q}dy.  \tag{5.2}
\end{equation}%
Since $A$ is symmetric and positive definite, we have%
\begin{equation}
\sup\limits_{k,i}\sum\limits_{j=1}^{N}\left\vert \frac{\lambda _{k}}{D\left(
\lambda _{k}\right) }A_{ji}\left( \lambda _{k}\right) \right\vert ^{q}\leq C.
\tag{5.3}
\end{equation}%
From $\left( 5.2\right) $ and $\left( 5.3\right) $ we get 
\begin{equation*}
\int\limits_{\Omega }\left\Vert \sum\limits_{k=1}^{\mu }r_{k}\left( y\right)
B\left( \lambda _{k}\right) u_{k}\right\Vert _{l_{q}}^{q}dy\leq
C\int\limits_{\Omega }\left\Vert \sum\limits_{k=1}^{\mu }r_{k}\left(
y\right) u_{k}\right\Vert _{l_{q}}^{q}dy.
\end{equation*}%
i.e., the operator $A$ is $R$-positive in $l_{q}.$ From Theorem 2.1 we
obtain that problem $\left( 5.1\right) $ has a unique solution $u\in $ $%
H_{p}^{m}\left( R^{n};l_{q}\left( A\right) ,l_{q}\right) $ for $f\in
L_{p}\left( R^{n};l_{q}\right) $ and the following estimate holds 
\begin{equation*}
\sum\limits_{\left\vert \alpha \right\vert \leq m}t\left( \alpha \right)
\left\vert \lambda \right\vert ^{1-\frac{\left\vert \alpha \right\vert }{m}%
}\left\Vert D^{\alpha }u\right\Vert _{L_{p}\left( R^{n};l_{q}\right)
}+\left\Vert Au\right\Vert _{L_{p}\left( R^{n};l_{q}\right) }\leq
M\left\Vert f\right\Vert _{L_{p}\left( R^{n};l_{q}\right) }.
\end{equation*}

From the above estimate we obtain the assertion.

\begin{center}
\textbf{Acknowledgements}
\end{center}

The author would like to express a gratitude to Dr. Erchan Aptoula for his
useful advices in English in preparing of this paper.

\vspace{3mm}

\vspace{3mm}

\textbf{References}

\ \ \ \ \ \ \ \ \ \ \ \ \ \ \ \ \ \ \ \ \ \ \ \ \ \ \ \ \ \ \ \ \ \ \ \ \ \
\ \ \ \ \ \ \ \ \ \ \ \ \ \ \ \ \ \ \ \ \ \ \ \ \ \ \ \ \ \ \ \ \ \ \ \ \ \
\ \ \ \ \ \ \ \ \ 

\begin{enumerate}
\item Amann H., Linear and quasi-linear equations,1, Birkhauser, Basel 1995.

\item Besov, O. V., P. Ilin, V. P., Nikolskii, S. M., Integral
representations of functions and embedding theorems, Nauka, Moscow, 1975.

\item Bourgain J., Some remarks on Banach spaces in which martingale
differences are unconditional, Arkiv Math. 21 (1983), 163-168.

\item Burkholder D. L., A geometrical conditions that implies the existence
certain singular integral of Banach space-valued functions, Proc. conf.
Harmonic analysis in honor of Antonu Zigmund, Chicago, 1981,Wads Worth,
Belmont, (1983), 270-286.

\item Dore C. and Yakubov S., Semigroup estimates and non coercive boundary
value problems, Semigroup Forum 60 (2000), 93-121.

\item Dubinskii, Yu. A., The Cauchy problem and pseudodifferential operators
in the complex domain, Uspekhi Mat. Nauk, 45:2(272) (1990), 115--142.

\item Diestel J., Jarchow H., Tonge A., Absolutely summing operators,
Cambridge Univ. Press, Cambridge, 1995.

\item Denk, R., Krainer T., $R$-Boundedness, pseudodifferential operators,
and maximal regularity for some classes of partial differential operators,
Manuscripta Math. 124(3) (2007),, 319-342.

\item Haller R., Heck H., Noll A., Mikhlin's theorem for operator-valued
Fourier multipliers in $n$ variables, Math. Nachr. 244 (2002), 110-130.

\item Hermander L., The analysis of linear partial differential operators,
V. 3, 4, Springer-Verlag, New York, 1985.

\item Favini A., Shakhmurov V. B, Yakov Yakubov, Regular boundary value
problems for complete second order elliptic differential-operator equations
in UMD Banach spaces, Semigroup form, 2009, v.79 (1), 22-54.

\item Krein S. G., \textquotedblright Linear differential equations in
Banach space\textquotedblright , American Mathematical Society, Providence,
1971.

\item Lizorkin P. I., Generalized Liouville differentiation and functional
spaces $L_{p}^{r}\left( E_{n}\right) .$ Embedding theorems, Mathematics of
the USSR-Sbornik (3)60 (1963), 325-353.

\item Nagase, M., The $L_{p}$-boundedness of pseudodifferential operators
with non-regular symbols, Comm. in Partial Differential Equations, 2 (10),
1045-1061 (1977).

\item Portal, P., \v{S}trkalj, \v{Z}., Pseudodifferential operators on
Bochner spaces and an application, Math. Z. (2006) 253: 805-819.

\item Shakhmurov V. B., Coercive boundary value problems for regular
degenerate differential-operator equations, J. Math. Anal. Appl., 292 ( 2)
(2004), 605-620.

\item Shakhmurov V. B., Embedding and separable differential operators in
Sobolev-Lions type spaces, Mathematical Notes, 2008, v. 84, no 6, 906-926.

\item Shakhmurov V. B., Maximal regular abstract elliptic equations and
applications, Siberian Mathematical Journal, v.51, no 5, 935-948, 2010.

\item Shakhmurov V. B., Embedding and maximal regular differential operators
in Banach-valued weighted spaces\textbf{, }Acta mathematica Sinica, 22(5)
(2006), 1493-1508.

\item Triebel H., \textquotedblright Interpolation theory, Function spaces,
Differential operators.\textquotedblright , North-Holland, Amsterdam, 1978.

\item Taylor, M., Pseudo-differential operators, Princeton University press,
Princeton, NJ, 1981.

\item Treves, F, Introduction to pseudo-differential and Fourier integral
operators, Plenum press, New York and London, 1980.

\item Weis L, Operator-valued Fourier multiplier theorems and maximal $L_{p}$
regularity, Math. Ann. 319, (2001), 735-758.

\item Yakubov S. and Yakubov Ya., \textquotedblright Differential-operator
Equations. Ordinary and Partial \ Differential Equations \textquotedblright
, Chapman and Hall /CRC, Boca Raton, 2000.
\end{enumerate}

\begin{quote}
\ 

\bigskip
\end{quote}

\end{document}